\documentstyle[twoside,amstex]{article}
\input amssym.def
\input amssym.tex
\def\bc{\begin{center}}
\def\ec{\end{center}}
\def\no{\noindent}

\begin{document}
\thispagestyle{empty} \vspace*{3 true cm} \pagestyle{myheadings}
\markboth {\hfill {\sl Huanyin Chen, O. Gurgun and H. Kose}\hfill}
{\hfill{\sl Uniquely Strongly Clean Triangular Matrix
Rings}\hfill} \vspace*{-1.5 true cm} \bc{\large\bf Uniquely
Strongly Clean Triangular Matrix Rings}\ec

\vskip6mm
\bc{{\bf Huanyin Chen}\\[1mm]
Department of Mathematics, Hangzhou Normal University\\
Hangzhou 310036, China, huanyinchen@@aliyun.com}\ec
\bc{{\bf O. Gurgun}\\[1mm]
Department of Mathematics, Ankara University\\ 06100 Ankara,
Turkey, orhangurgun@@gmail.com}\ec
\bc{{\bf H. Kose}\\[1mm]
Department of Mathematics, Ahi Evran University\\ Kirsehir,
Turkey, handankose@@gmail.com}\ec

\begin{abstract} A ring $R$ is uniquely (strongly) clean provided that for any $a\in R$ there exists a unique idempotent $e\in R \big( \in comm(a)\big)$ such that
$a-e\in U(R)$. Let $R$ be a uniquely bleached ring. We prove, in
this note, that $R$ is uniquely clean if and only if $R$ is
abelian, and $T_n(R)$ is uniquely strongly clean for all $n\geq
1$, if and only if $R$ is abelian, $T_n(R)$ is uniquely strongly
clean for some $n\geq 1$. In the commutative case, the more
explicit results are obtained. These also generalize the main
theorems in [6] and [7], and
provide many new class of such rings.\\[2mm]
\end{abstract}

\section{Introduction}

The commutant of an element $a$ in a ring $R$ is defined by
$comm(a)=\{x\in R~|~xa=ax\}$. A ring $R$ is strongly clean
provided that for any $a\in R$ there exists an idempotent $e\in
comm(a)$ such that $a-e\in U(R)$. Strongly clean triangular matrix
rings are extensively studied by many authors, e.g., [1] and [3].
A ring $R$ is called uniquely clean provided that for any $a\in R$
there exists a unique idempotent $e\in R$ such that $a-e\in U(R)$.
Many characterizations of such rings are studied in [2-4] and [8].
Following J. Chen; Z. Wang and Y. Zhou [6], a ring $R$ is called
uniquely strongly clean provided that for any $a\in R$ there
exists a unique idempotent $e\in comm(a)$ such that $a-e\in U(R)$.
Uniquely strongly cleanness behaves very different from the
properties of uniquely clean rings (cf. [6]). In general, matrix
rings have not such properties (cf. [10, Proposition 11.8]). Thus,
it is attractive to investigate uniquely strongly cleanness of
triangular matrix rings. Chen et al. proved that if $R$ is
commutative then $R$ is uniquely clean if and only if $T_n(R)$ is
uniquely strongly clean for all $n\geq 1$, if and only if $T_n(R)$
is uniquely strongly clean for some $n\geq 1$.

[6, Question 12] and [10, Question 11.13] asked that if
"commutative" in the preceding result can be replaced by
"abelian". The motivation of this note article is to explore this
problem. Recall that a ring $R$ is uniquely bleached provided that
for any $a\in J(R),b\in U(R)$, $l_a-r_b$ and $l_b-r_a$ are
isomorphism. Let $R$ be a uniquely bleached ring. We prove, in
this note, that $R$ is uniquely clean if and only if $R$ is
abelian, and $T_n(R)$ is uniquely strongly clean for all $n\geq
1$, if and only if $R$ is abelian, $T_n(R)$ is uniquely strongly
clean for some $n\geq 1$. In the commutative case, the more
explicit results are obtained. These also generalize the main
theorems in [6] and [7], and provide many new class of such rings.

We write $U(R)$ for the set of all invertible elements in $R$ and
$J(R)$ for the Jacobson radical of $R$. $T_n(R)$ stand for the
rings of all $n\times n$ triangular matrices over a ring $R$. Let
$a,b\in R$. We denote the map from $R$ to $R: x\mapsto ax-xb$ by
$l_a-r_b$.

\section{The Main Results}

\vskip4mm Clearly,we have $\{$ uniquely clean rings $\}\subsetneq
\{$ abelian clean rings $\}\subsetneq \{$ uniquely strongly clean
rings $\}$. We begin with

\vskip4mm \hspace{-1.8em} {\bf Theorem 1.}\ \ {\it Let $R$ be a
uniquely bleached ring. Then the following are
equivalent:}\vspace{-.5mm}
\begin{enumerate}
\item [(1)]{\it $R$ is uniquely clean.}
\vspace{-.5mm}
\item [(2)]{\it $R$ is abelian, and $T_n(R)$ is uniquely strongly clean for all $n\geq 1$.}
\vspace{-.5mm}
\item [(3)]{\it $R$ is abelian, $T_n(R)$ is uniquely strongly clean for some $n\geq 1$.}
\end{enumerate}\vspace{-.5mm} {\it Proof.}\ \ $(1)\Rightarrow (2)$
In view of [8, Theorem 20], $R$ is abelian.  Clearly, the result
holds for $n=1$. Assume that the result holds for $n (n\geq 1)$.
Let $A=\left(
\begin{array}{cc}
a_{11}&\alpha\\
&A_{1}
\end{array}
\right)\in T_{n+1}(R)$ where $a_{11}\in R, \alpha\in M_{1\times
n}(R)$ and $A_1\in T_{n}(R)$. By hypothesis, we can find a unique
idempotent $e_{11}\in R$ such that $u_{11}:=a_{11}-e_{11}\in U(R)$
and $a_{11}e_{11}=e_{11}a_{11}$. Furthermore, we have a unique
idempotent $E_1\in T_{n}(R)$ such that $U_1:=A_1-E_1\in
U\big(T_{n}(R)\big)$ and $A_1E_1=E_1A_1$; hence, $U_1E_1=E_1U_1$.
Let $E=\left(
\begin{array}{cc}
e_{11}&x\\
&E_1
\end{array}
\right)$ and $U=\left(
\begin{array}{cc}
u_{11}&\alpha-x\\
&U_1
\end{array}
\right)$, where $x\in M_{1\times n}(R)$. Observing that
$$\begin{array}{c}
E^2=E\Leftrightarrow e_{11}x+xE_1=x;\\
UE=EU\Leftrightarrow u_{11}x+(\alpha-x)E_1=e_{11}(\alpha-x)+xU_1,
\end{array}$$
and then $$(u_{11}+2e_{11}-1)x-xU_1=e_{11}\alpha-\alpha
E_1~~(*).$$ In view of [8, Theorem 20], $R/J(R)$ is Boolean, and
so $2\in J(R)$. Furthermore, $u_{11}\in 1+J(R)$. This shows that
$u_{11}+2e_{11}-1\in J(R)$. Write $$\begin{array}{c} x=\left(
\begin{array}{ccc}
x_1&\cdots &x_n
\end{array}
\right), e_{11}\alpha-\alpha E_1=\left(
\begin{array}{ccc}
v_1&\cdots &v_n
\end{array}
\right),\\
U_1=\left(
\begin{array}{cccc}
c_{11}&c_{12}&\cdots &c_{1n}\\
&c_{22}&\cdots&c_{2n}\\
&&\ddots&\\
&&&c_{nn}
\end{array}
\right), \mbox{where each}~c_{ii}\in U(R). \end{array}$$ Then
$$\begin{array}{c}
(u_{11}+2e_{11}-1)x_1-x_1c_{11}=v_1;\\
(u_{11}+2e_{11}-1)x_2-x_2c_{22}=v_1+x_1c_{12};\\
\vdots\\
(u_{11}+2e_{11}-1)x_n-x_nc_{nn}=v_n+x_1c_{1n}+\cdots
+x_{n-1}c_{(n-1)n}.
\end{array}$$ By hypothesis, we have a unique $x$ such that $(*)$
holds. As $E_1U_1=U_1E_1$, we get
$$\begin{array}{ll}
&(u_{11}+2e_{11}-1)x(e_{11}I_n+E_1)-x(e_{11}I_n+E_1)U_1\\
=&\alpha (e_{11}I_n-E_1)(e_{11}I_n+E_1)\\
=&\alpha (e_{11}I_n-E_1)\\
=&(u_{11}+2e_{11}-1)x-xU_1.
\end{array}$$ Set $y=x(e_{11}I_n+E_1)$. Then
$(u_{11}+2e_{11}-1)(y-x)-(y-x)U_1=0$.

Write $y-x=\left(
\begin{array}{ccc}
z_1&\cdots &z_n
\end{array}
\right)$. Then $$\begin{array}{c}
(u_{11}+2e_{11}-1)z_1-z_1c_{11}=0;\\
(u_{11}+2e_{11}-1)z_2-z_2c_{22}=z_1c_{12};\\
\vdots\\
(u_{11}+2e_{11}-1)z_n-z_nc_{nn}=z_1c_{1n}+\cdots
+z_{n-1}c_{(n-1)n}.
\end{array}$$ By hypothesis, we get each $z_i=0$, and so $y=x$.
This implies that $e_{11}x+xE_1=x$. Furthermore,
$u_{11}x-(\alpha-x)E_1=e_{11}(\alpha-x)+xU_1.$ Therefore, we have
a uniquely strongly clean expression
$$A=\left(
\begin{array}{cc}
e_{11}&x\\
&E_1
\end{array}
\right)+\left(
\begin{array}{cc}
a_{11}-e_{11}&\alpha-x\\
&A_1-E_1
\end{array}
\right).
$$ By induction, $T_n(R)$
is uniquely strongly clean for all $n\in {\Bbb N}$.

$(2)\Rightarrow (3)$ is trivial.

$(3)\Rightarrow (1)$ In light of [6, Example 5], $R$ is uniquely
strongly clean, hence the result by [8, Theorem 20].\hfill$\Box$

\vskip4mm \hspace{-1.8em} {\bf Corollary 2.}\ \ {\it Let $R$ be a
ring with nil Jacobson radical. Then the following are
equivalent:}\vspace{-.5mm}
\begin{enumerate}
\item [(1)]{\it $R$ is uniquely clean.}
\vspace{-.5mm}
\item [(2)]{\it $R$ is abelian, and $T_n(R)$ is uniquely strongly clean for all $n\geq 1$.}
\vspace{-.5mm}
\item [(3)]{\it $R$ is abelian, and $T_n(R)$ is uniquely strongly clean for some $n\geq 1$.}
\end{enumerate}\vspace{-.5mm} {\it Proof.}\ \ Let $a\in U(R), b\in J(R)$. Write $b^n=0$. Choose
$\varphi=l_{a^{-1}}+l_{a^{-2}}r_{b}+\cdots +l_{a^{-n}}r_{b^{n-1}}:
R\to R$. One easily checks that $\big(l_a-r_b\big)\varphi=\varphi
\big(l_a-r_b\big)=1_R$. Thus, $l_a-r_b: R\to R$ are isomorphism.
Similarly, $l_b-r_a$ is isomorphic. Therefore $R$ is uniquely
bleached, and the result follows by Theorem 1.\hfill$\Box$

\vskip4mm \hspace{-1.8em} {\bf Corollary 3.}\ \ {\it Let $R$ be a
ring for which some power of each element in $J(R)$ is central.
Then the following are equivalent:}\vspace{-.5mm}
\begin{enumerate}
\item [(1)]{\it $R$ is uniquely clean.}
\vspace{-.5mm}
\item [(2)]{\it $R$ is abelian, and $T_n(R)$ is uniquely strongly clean for all $n\geq 1$.}
\vspace{-.5mm}
\item [(3)]{\it $R$ is abelian, and $T_n(R)$ is uniquely strongly clean for some $n\geq 1$.}
\end{enumerate}\vspace{-.5mm} {\it Proof.}\ \ Let $a\in U(R), b\in J(R)$. Write $b^n\in C(R)$. Choose
$\varphi=l_{a^{-1}}+l_{a^{-2}}r_{b}+\cdots +l_{a^{-n}}r_{b^{n-1}}:
R\to R$. Then
$\big(l_a-r_b\big)\varphi=l_{1}-l_{a^{-n}}r_{b^n}=l_{1-a^{-n}b^n}$.
Further, we check that $\varphi\big(l_a-r_b\big)=l_{1-a^{-n}b^n}$.
Thus, $l_a-r_b: R\to R$ is isomorphic. Likewise, $l_b-r_a$ is
isomorphic. That is, $R$ is uniquely bleached. Therefore we
complete the proof by Theorem 1.\hfill$\Box$

\vskip4mm Immediately, we show that every triangular matrix ring
over a Boolean ring is uniquely strongly clean. For instance,
$T_n\big({\Bbb Z}_2\big) (n\geq 2)$ is uniquely strongly clean,
while it is not uniquely clean.

\vskip4mm \hspace{-1.8em} {\bf Corollary 4 [6, Theorem 10].}\ \
{\it Let $R$ be a commutative ring. Then the following are
equivalent:}\vspace{-.5mm}
\begin{enumerate}
\item [(1)]{\it $R$ is uniquely clean.}
\vspace{-.5mm}
\item [(2)]{\it $T_n(R)$ is uniquely strongly clean for all $n\geq 1$.}
\vspace{-.5mm}
\item [(3)]{\it $T_n(R)$ is uniquely strongly clean for some $n\geq 1$.}
\end{enumerate}\vspace{-.5mm} {\it Proof.}\ \ It is clear as every
commutative ring is uniquely bleached and abelian.\hfill$\Box$

\vskip4mm \hspace{-1.8em} {\bf Lemma 5.}\ \ {\it Let $R$ be a
ring. If $T_2(R)$ is uniquely strongly clean, then $R$ is uniquely
bleached.}\vskip2mm\hspace{-1.8em} {\it Proof.}\ \ In view of [6,
Example 5], $R$ is uniquely strongly clean. Let $a\in J(R)$ and
$b\in U(R)$, and let $r\in R$. Choose $A=\left(
\begin{array}{cc}
a&-r\\
&b
\end{array}
\right)\in T_2(R)$. Then there exists a unique idempotent
$E=(e_{ij})\in T_2(R)$ such that $A-E\in U\big(T_2(R)\big)$ and
$EA=AE$. Clearly, $e_{11},e_{22}\in R$ are idempotents. Further,
$a-e_{11}\in U(R)$ and $b-e_{22}\in U(R)$. As $a-0\in U(R)$ and
$b-1\in U(R)$, by the unique strong cleanness of $R$, we get
$e_{11}=0$ and $e_{22}=1$. Thus, $E=\left(
\begin{array}{cc}
0&x\\
&1
\end{array}
\right)$ for some $x\in R$. It follows from $EA=AE$ that
$ax-xb=r$. Assume that $ay-yb=r$. Then we have an idempotent
$F=\left(
\begin{array}{cc}
0&y\\
&1
\end{array}
\right)$ such that $A-F\in U\big(T_2(R)\big)$ and $AF=FA$. By the
uniqueness of $E$, we get $x=y$. Therefore, $l_a-r_b: R\to R$ is
an isomorphism. Likewise, $l_b-r_a: R\to R$ is an isomorphism.
Accordingly, $R$ is uniquely bleached, as asserted.\hfill$\Box$

\vskip4mm \hspace{-1.8em} {\bf Theorem 6.}\ \ {\it Let $R$ be an
abelian ring. Then the following are equivalent:}\vspace{-.5mm}
\begin{enumerate}
\item [(1)]{\it $R$ is uniquely bleached, uniquely clean.}
\vspace{-.5mm}
\item [(2)]{\it $T_n(R)$ is uniquely strongly clean for all $n\geq 1$.}
\vspace{-.5mm}
\item [(3)]{\it $T_n(R)$ is uniquely strongly clean for some $n\geq 2$.}
\end{enumerate}\vspace{-.5mm} {\it Proof.}\ \ $(1)\Rightarrow (2)$
is proved by Theorem 1.

$(2)\Rightarrow (3)$ is trivial.

$(3)\Rightarrow (1)$ In view of [6, Example 5], $T_2(R)$ is
uniquely strongly clean. Further, $R$ is uniquely strongly clean.
As $R$ is abelian, $R$ is uniquely clean. Therefore the proof is
complete by Lemma 5.\hfill$\Box$

\vskip4mm Let $Lat(R)$ be the lattice\index{distribute lattice} of
all right ideals of a ring $R$. A ring $R$ said to be a $D$-ring
if $Lat(R)$ is distribute, i.e., $A\bigcap (B+C)=(A\bigcap
B)+(A\bigcap C)$ for all $A,B,C\in Lat(R)$.

\vskip4mm \hspace{-1.8em} {\bf Corollary 7.}\ \ {\it Let $R$ be a
$D$-ring. Then the following are equivalent:}\vspace{-.5mm}
\begin{enumerate}
\item [(1)]{\it $R$ is uniquely bleached, uniquely clean.}
\vspace{-.5mm}
\item [(2)]{\it $T_n(R)$ is uniquely strongly clean for all $n\geq 1$.}
\vspace{-.5mm}
\item [(3)]{\it $T_n(R)$ is uniquely strongly clean for some $n\geq 1$.}
\end{enumerate}\vspace{-.5mm} {\it Proof.}\ \ Let
$e\in R$ be an idempotent, and let $a\in R$. Choose $f=e-(1-e)ae$.
Clearly, $eR,fR\in Lat(R)$. Then $eR=\big(eR\bigcap
fR\big)+\big(eR\bigcap (1-f)R\big)$; hence, $feR\subseteq eR$.
Thus, $efe=fe=f$. This implies that $ae=eae$. Similarly, we show
that $ea=eae$. Therefore $R$ is abelian, and the proof is complete
by Theorem 6.\hfill$\Box$

\vskip4mm \hspace{-1.8em} {\bf Corollary 8 ([7, Theorem 1]).}\ \
{\it Let $R$ be a local ring. Then the following are
equivalent:}\vspace{-.5mm}
\begin{enumerate}
\item [(1)]{\it $R$ is uniquely bleached and $R/J(R)\cong {\Bbb Z}_2$.}
\vspace{-.5mm}
\item [(2)]{\it $T_n(R)$ is uniquely strongly clean for all $n\geq 1$.}
\vspace{-.5mm}
\item [(3)]{\it $T_n(R)$ is uniquely strongly clean for some $n\geq 2$.}
\end{enumerate}\vspace{-.5mm} {\it Proof.}\ \ As $R$ is local, $R$ is uniquely clean if and only if $R/J(R)\cong {\Bbb Z}_2$, by [8, Theorem 15]. Therefore we complete the proof
from Theorem 6.\hfill$\Box$

\vskip4mm The double commutant of an element $a$ in a ring $R$ is
defined by $comm^2(a)=\{x\in R~|~xy=yx~\mbox{for all}~y\in
comm(a)\}$. Clearly, $comm^2(a)\subseteq comm(a)$. We end this
note by a more explicit result than [6, Theorem 10].

\vskip4mm \hspace{-1.8em} {\bf Theorem 9.}\ \ {\it Let $R$ be a
commutative ring, and let $n\in {\Bbb N}$. Then the following are
equivalent:}\vspace{-.5mm}
\begin{enumerate}
\item [(1)]{\it $R$ is uniquely clean.}
\vspace{-.5mm}
\item [(2)]{\it For any $A\in T_n(R)$, there exists a unique idempotent $E\in comm^2(A)$ such that
$A-E\in U\big(T_n(R)\big)$.}\end{enumerate}\vspace{-.5mm} {\it
Proof.}\ \ $(1)\Rightarrow (2)$ For any $A\in T_n(R)$, we claim
that there exists an idempotent $E\in comm^2(A)$ such that $A-E\in
J\big(T_2(R)\big)$. As $R$ is a commutative uniquely clean ring,
the result holds for $n=1$, by [8, Theorem 20]. Suppose that the
result holds for $n-1 (n\geq 2)$. Let $A=\left(
\begin{array}{cc}
a_{11}&\alpha\\
&A_{1}
\end{array}
\right)\in T_n(R)$ where $a_{11}\in R, \alpha\in M_{1\times
(n-1)}(R)$ and $A_1\in T_{n-1}(R)$. Then there is an idempotent
$e_{11}\in R$ such that $w_{11}:=a_{11}-e_{11}\in J(R)$. By
hypothesis, there is an idempotent $E_1\in T_{n-1}(R)$ such that
$W_1:=A_1-E_1\in J\big(T_{n-1}(R)\big)$ and $E_1\in comm^2(A_1)$.
In view of [8, Theorem 20], $R/J(R)$ is Boolean, and so $2\in
J(R)$. Thus, $W_1+\big(1-2e_{11}-w_{11}\big)I_{n-1}\in
I_{n-1}+J\big(T_{n-1}(R)\big)\subseteq U\big(T_{n-1}(R)\big)$. Let
$E=\left(
\begin{array}{cc}
e_{11}&\beta\\
&E_1
\end{array}
\right)$, where $\beta=\alpha (E_1-e_{11}I_{n-1})
\big(W_1+(1-2e_{11}-w_{11})I_{n-1}\big)^{-1}$. Then $A-E\in
J\big(T_n(R)\big)$. As $e_{11}\beta+\beta E_1=\beta
(E_1+e_{11}I_{n-1})=\alpha (E_1-e_{11}I_{n-1})(E_1+e_{11}I_{n-1})
\big(W_1+(1-2e_{11}-w_{11})I_{n-1}\big)^{-1}=\beta$, we see that
$E=E^2$.

For any $X=\left(
\begin{array}{cc}
x_{11}&\gamma\\
&X_1
\end{array}
\right)\in comm(A)$, we have $\alpha (X_1-x_{11}I_{n-1})=\gamma
(A_1-a_{11}I_{n-1}).$ Since $E_1\in comm^2(A_1)$, as in the proof
of [4, Theorem 4.11], we check that
$$\begin{array}{ll}
&\gamma(A_1-a_{11}I_{n-1})(E_{1}-e_{11}I_{n-1})\\
=&\alpha (X_1-x_{11}I_{n-1})(E_{1}-e_{11}I_{n-1})\\
=&\alpha (E_{1}-e_{11}I_{n-1})(X_1-x_{11}I_{n-1})\\
=&\beta \big(W_1+(1-2e_{11}-w_{11})I_{n-1}\big)(X_1-x_{11}I_{n-1})\\
=&\beta
(X_1-x_{11}I_{n-1})\big(W_1+(1-2e_{11}-w_{11})I_{n-1}\big).
\end{array}$$
Moreover,
$$\begin{array}{ll}
&\gamma(A_1-a_{11}I_{n-1})(E_{1}-e_{11}I_{n-1})\\
=&\gamma(E_{1}-e_{11}I_{n-1})\big(E_1+W_1-(e_{11}+w_{11})I_{n-1}\big)\\
=&\gamma(E_{1}-e_{11}I_{n-1})\big(E_1+e_{11}I_{n-1}+(W_1-2e_{11}-w_{11})I_{n-1}\big)\\
=&\gamma\big(E_1-e_{11}I_{n-1}+(E_{1}-e_{11}I_{n-1})(W_1-2e_{11}-w_{11})I_{n-1}\big)\\
=&\gamma
(E_1-e_{11}I_{n-1})\big(W_1+(1-2e_{11}-w_{11})I_{n-1}\big).
\end{array}$$
This shows that $\gamma (E_1-e_{11}I_{n-1})=\beta
(X_1-x_{11}I_{n-1})$. Thus, we get $e_{11}\gamma+\beta
X_{1}=x_{11}\beta+\gamma E_1$; hence, $EX=XE$. That is, $E\in
comm^2(A)$, as claimed.

Therefore, for any $A\in T_n(R)$, there exists an idempotent $E\in
comm^2(A)$ such that $W:=A-E\in J\big(T_2(R)\big)$. Hence,
$A=\big(I_n-E\big)+\big(2E-I_n+W\big)$. As
$\big(2E-I_n\big)^2=I_n$, we see that $2E-I_n\in
U\big(T_n(R)\big)$. Additionally, $I_n-E\in comm^2(A)$. Suppose
that there exists an idempotent $F\in comm^2(A)$ such that $A-F\in
U\big(T_2(R)\big)$. Then $I_n-E,F\in comm(A)$. In view of
Corollary 4, $T_n(R)$ is uniquely strongly clean. This implies
that $F=I_n-E$, proving $(1)$.

$(2)\Rightarrow (1)$ Let $a\in R$. Then $A=diag(a,a,\cdots ,a)\in
T_n(R)$. Then we can find a unique idempotent $E=(e_{ij})\in
comm^2(A)$ such that $A-E\in U\big(T_n(R)\big)$. This implies that
$e_{11}\in R$ is an idempotent and $a-e_{11}\in U(R)$. Suppose
that $a-e\in U(R)$ with an idempotent $e\in R$. Then
$F=diag(e,e,\cdots ,e)\in T_n(R)$ is an idempotent. Further, $F\in
comm^2(A)$, and that $A-F\in U\big(T_n(R)\big)$. By the
uniqueness, we get $E=F$, and then $e=e_{11}$. Therefore $R$ is
uniquely clean, as asserted.\hfill$\Box$

\vskip4mm \hspace{-1.8em} {\bf Corollary 10.}\ \ {\it Let $R$ be a
commutative uniquely clean ring, and let $A\in T_n(R)$. Then for
any $A\in T_n(R)$, there exists a unique idempotent $E\in comm(A)$
such that $A-E, A+E\in U\big(T_n(R)\big)$.}\vskip2mm
\hspace{-1.8em} {\it Proof.}\ \ In view of Theorem 9, we have a
unique idempotent $E\in comm^2(A^2)$ such that $U:=A^2-E\in
GL_n(R)$. Explicit, we may assume that $U\in
-I_n+J\big(T_n(R)\big)$. Thus, $(A-E)(A+E)\in U\big(T_n(R)\big)$.
Therefore, $A-E, A+E\in U\big(T_n(R)\big)$. Assume that there is
an idempotent $F\in comm(A)$ such that $A-F,A+F\in
U\big(T_n(R)\big)$. Then $A^2-F=(A-F)(A+F)\in U\big(T_n(R)\big)$.
Hence, $(I_n-E)+F=(A^2+I_n-E)-(A^2-F)\in U\big(T_n(R)\big)$.
Clearly, $(I_n-E)F=F(I_n-E)$, and so
$\big((I_n-E)-F\big)\big((I_n-E)+F-I_n\big)=0$. This shows that
$F=I_n-E$, as desired.\hfill$\Box$

\vskip4mm \hspace{-1.8em} {\bf Corollary 11.}\ \ {\it Let $R$ be a
Boolean ring, and let $A\in T_n(R)$.}\vspace{-.5mm}
\begin{enumerate}
\item [(1)]{\it There exists a unique
idempotent $E\in comm^2(A)$ such that $A-E\in U\big(T_n(R)\big)$.}
\vspace{-.5mm}
\item [(2)]{\it There exists a unique
idempotent $E\in comm^2(A)$ such that $A-E, A+E\in
U\big(T_n(R)\big)$.}\end{enumerate}\vspace{-.5mm} {\it Proof.}\ \
Clearly, $R$ is a commutative uniquely clean ring, hence the
result by Theorem 9 and Corollary 10.\hfill$\Box$

\vskip15mm\no {\Large\bf Acknowledgements} \vskip4mm This research
was supported by the Natural Science Foundation of Zhejiang
Province (LY13A0 10019) and the Scientific and Technological
Research Council of Turkey (2221 Visiting Scientists Fellowship
Programme).

\end{document}